\documentclass{amsart}

\usepackage{amssymb}
\usepackage[all]{xy}

\newtheorem{thm}{Theorem}

\newtheorem{lem}[thm]{Lemma}
\newtheorem{cor}[thm]{Corollary}

\newtheorem{prop}[thm]{Proposition}

\theoremstyle{definition}
\newtheorem{defn}[thm]{Definition}

\newtheorem{say}[thm]{}
\newtheorem{exmp}[thm]{Example}

\newtheorem{ques}[thm]{Question}    

\newtheorem{rem}[thm]{Remark}          

\newtheorem*{ack}{Acknowledgments}      
\newtheorem{notation}[thm]{Notation}   
  
\newtheorem{defn-thm}[thm]{Definition--Theorem}  
\newtheorem{defn-lem}[thm]{Definition--Lemma}  

\theoremstyle{remark}

\renewcommand{\c}[0]{{\mathbb C}}  

\renewcommand{\o}[0]{{\mathcal O}} 
\newcommand{\z}[0]{{\mathbb Z}}

\renewcommand{\r}[0]{{\mathbb R}} 

\renewcommand{\a}[0]{{\mathbb A}}

\newcommand{\s}[0]{{\mathbb S}}

\newcommand{\p}[0]{{\mathbb P}}

\newcommand{\map}[0]{\dasharrow}
\newcommand{\qtq}[1]{\quad\mbox{#1}\quad}
\newcommand{\spec}[0]{\operatorname{Spec}}
\newcommand{\pic}[0]{\operatorname{Pic}}

\newcommand{\gal}[0]{\operatorname{Gal}}

\newcommand{\rank}[0]{\operatorname{rank}}

\newcommand{\aut}[0]{\operatorname{Aut}}

\newcommand{\sing}[0]{\operatorname{Sing}}

\newcommand{\disc}[0]{\operatorname{discr}} 
 
\newcommand{\chow}[0]{\operatorname{Chow}}
\newcommand{\chr}[0]{\operatorname{char}}

\newcommand{\hilb}[0]{\operatorname{Hilb}}
\newcommand{\grass}[0]{\operatorname{Grass}}

\newcommand{\univ}[0]{\operatorname{Univ}}

\newcommand{\rdown}[1]{\lfloor{#1}\rfloor}

\newcommand{\onto}[0]{\twoheadrightarrow}

\newcommand{\ord}[0]{\operatorname{ord}}

\newcommand{\isom}[0]{\operatorname{Isom}}

\newcommand{\mor}[0]{\operatorname{Mor}}

\newcommand{\ind}[0]{\operatorname{index}}

\newcommand{\rc}[0]{\operatorname{RatCurve}}

\def\into{\DOTSB\lhook\joinrel\to}

\def\loccoh#1.#2.#3.#4.{H^{#1}_{#2}(#3,#4)}

\DeclareMathAlphabet{\mathchanc}{OT1}{pzc}%
                                {m}{it}

\newcommand{\simb}[0]{\stackrel{bir}{\sim}}
\newcommand{\sims}[0]{\stackrel{stab}{\sim}}

\newcommand{\PGL}{\mathrm{PGL}}

\newcommand{\OO}{\mathrm{O}}
\newcommand{\PGO}{\mathrm{PGO}}

\newcommand{\OG}{\mathrm{OG}}
\newcommand{\OF}{\mathrm{OF}}
\newcommand{\G}{\mathrm{G}}

\newcommand{\PTC}{\mathrm{PTC}}

\newcommand{\sym}[0]{\operatorname{Sym}}

\usepackage[all]{xy}\xyoption{dvips}

\newcommand{\norm}[0]{\operatorname{norm}}

\newcommand{\emb}[0]{\operatorname{Emb}}

\newcommand{\maps}[0]{\operatorname{Map}}

\newcommand{\mo}[0]{\operatorname{M}^{\rm \circ}}

\begin{document}
\bibliographystyle{amsalpha}

{\ }

 \title[Symmetric powers of Severi-Brauer varieties]{Symmetric powers of Severi-Brauer varieties}
 \author{J\'anos Koll\'ar}
\begin{abstract} We classify products of symmetric  powers of a Severi-Brauer variety, up-to stable birational equivalence. The description also includes Grassmannians and moduli spaces of genus 0 stable maps.
\end{abstract}

 \maketitle

There are several ways to associate other varieties to
a Severi-Brauer variety $P$ defined over a field $k$.  These include
\begin{itemize}
\item the Grassmannians  $\grass(\p^{m-1}, P)$,
\item  the symmetric powers  $\sym^m(P)$ and
\item the moduli spaces $\bar M_0(P,d)$ of genus 0 stable maps of degree $d$ to $P$.
\end{itemize}
While all  these varieties are geometrically rational, they are usually not rational over the ground field and it is an interesting problem to understand
their birational properties over $k$. 
The results of this note  are partly weaker---since we describe only the 
stable birational equivalence classes---but partly stronger---since 
we also describe the products of these varieties.

Thus let  $\operatorname{MSym}(P)$ denote the  multiplicative 
monoid generated by 
 stable birational equivalence classes of  Grassmannians of $P$,
symmetric powers  of $P$ and  the moduli spaces   $\bar M_0(P,d)$.
We show that $\operatorname{MSym}(P)$ is finite,
identify its elements and also the multiplication rules. 

 Birational equivalence of two varieties is denoted by $X\simb Y$ and 
 stable birational equivalence by  $X\sims Y$.
See  Paragraph \ref{SB.basic.say} for the definition and basic properties of 
Severi-Brauer varieties.


\begin{thm}\label{sym.SB.thm}
 Let $P$ be a Severi-Brauer variety of  index $i=i(P)$. 
Then
\begin{enumerate}
\item $\operatorname{MSym}(P)=\{\grass(\p^{d-1},P): d\mid i(P)\}$ and
products are given by
\item $\grass(\p^{d-1},P)\times \grass(\p^{e-1},P)\sims \grass(\p^{(d,e)-1},P)$,
where  $(d,e)$  denotes the greatest common divisor.
The identity is  $\grass(\p^{i-1},P)\sims \p^0$.
\end{enumerate}
The class of an arbitrary Grassmannian is given by the rule
\begin{enumerate}\setcounter{enumi}{2}
\item $\grass(\p^{d-1},P)\sims  \grass(\p^{(d, i)-1},P)$.
\end{enumerate}
The class of an arbitrary symmetric power is given by the rules
\begin{enumerate}\setcounter{enumi}{3}
\item $\sym^d(P)\sims \sym^{(d, i)}(P)$ for every $d\geq 0$,
\item $\sym^d(P)\simb \grass(\p^{d-1}, P)\times \p^{d(d-1)}$ for $d\leq n+1$ and
\item $\sym^d(P)\times \sym^e(P)\sims \sym^{(d,e)}(P)$.
\end{enumerate}
The class of  $\bar M_0( P,d)$  is determined by the parity of $d$:
\begin{enumerate}\setcounter{enumi}{6}
\item $\bar M_0(P, 2e)\sims P$  save when $\dim P=e=1$ and 
\item $\bar M_0(P, 2e+1)\sims \grass(\p^1, P)$. Note that 
 $\grass(\p^1, P)$  is rational iff $i(P)\in\{1, 2\}$
and stably birational to $P$ iff $i(P)$ is odd.
\end{enumerate}
\end{thm}

The most natural description seems to be in terms of symmetric powers, so we start with them. The relationship with Grassmannians is easy to establish.
The moduli spaces $\bar M_0(P, d)$ end up birationally the simplest but
understanding them  is more subtle.  

The case $\dim P=e=1$  is exceptional in (\ref{sym.SB.thm}.7).
$\bar M_0(P, 2)$ aims to classify double covers of $P$ ramified at 2 points. 
The coarse moduli space is $\sym^2(P)\cong \p^2$. However, 
if $P\ncong \p^1$ then there are no such double covers defined over $k$.
The problem is that every double cover has an order 2 automorphism.
In all other cases, a dense open subset of   $\bar M_0(P, d)$
parametrizes maps without automorphisms, even embeddings if $\dim P\geq 3$.

\begin{say}[Severi-Brauer varieties I]\label{SB.basic.say}
Let $k$ be a field with separable closure $k^s$. A $k$-scheme $P$ is called 
a {\it Severi-Brauer variety}  if $P_{k^s}:=P\times_{\spec k}\spec k^s \cong \p^n$
for some $n$. We say that $P$ is trivial if $P\cong \p^n$. 
The {\it index} of $P$ is the gcd of all 0-cycles on $P$;
 it is denoted by $i(P)$.
A subscheme $L\subset P$ is called {\it twisted linear}
if $L_{k^s}$ is a linear subspace of $P_{k^s}\cong\p^n$.
Thus $L$ is also a Severi-Brauer variety. 
For a reduced subscheme $Y\subset X$ let $\langle Y\rangle$ denote the
minimal twisted linear subvariety containing $Y$. Thus
$\langle Y\rangle_{k^s}$ is the linear span of
$Y_{k^s}$.

The following basic results go back to Severi and Ch\^atelet, see  
\cite[Chap.5]{MR2266528} for a modern treatment and references.
\begin{enumerate}
\item $P$ is trivial iff $P(k)\neq \emptyset$.
\item  $i(P)$ divides $\dim P+1$,
\item $P$ has a 0-cycle  $Z$ of degree  $i(P)$ and 
$\dim \langle Z\rangle=i(P)-1$.
\item The minimal twisted linear subvarieties have dimension $i(P)-1$ and they are  isomorphic to each other; call this isomorphism class $P^{\rm min}$.
\item Given $P^{\rm min}$ and $r\geq 1$ there is a unique (up-to isomorphism)
Severi-Brauer variety $P_r$  of dimension  $r(\dim P^{\rm min}+1)-1$
such that $(P_r)^{\rm min}\cong P^{\rm min}$.
\item $P_r\simb P^{\rm min}\times \p^m$ for $m=(r-1)(\dim P^{\rm min}+1)$.
\item Two Severi-Brauer varieties $P_1, P_2$ are {\it Brauer-equivalent,}
denoted by $P_1\sim P_2$, iff
$P_1^{\rm min}\cong P_2^{\rm min}$. This holds iff the smaller dimensional one  is
isomorphic to a twisted linear subvariety of the other.
\end{enumerate}
\end{say}

\section{Symmetric powers}

A key step in understanding symmetric powers is the following.

\begin{thm} \cite{MR2090670} \label{kra-sal.thm}
Let $P$ be a Severi-Brauer variety of dimension $n$. 
Then $\sym^{n+1}(P)$ is rational.
\end{thm}

The following is a short geometric proof. The Euler number of $\p^n$ is $n+1$, thus a general
section of the tangent bundle $T_{\p^n}$ vanishes at $n+1$ points.
For any Severi-Brauer variety
this gives a dominant map  $\pi:H^0(P, T_P)\map \sym^{n+1}(P)$.

Let $Z\subset P$ be a reduced 0-cycle of degree $n+1$. 
Then  $\pi^{-1}(Z)$ is the linear space
$H^0\bigl(P, T_P(-Z)\bigr)\subset H^0(P, T_P)$ of dimension $n+1$.
Let $V\subset H^0(P, T_P)$ be a general affine-linear subspace of codimension
$n+1$. Then  $\pi|_V:V\map \sym^{n+1}(P)$ is birational. \qed
\medskip


\begin{cor} \label{kra-sal.cor}
Let $P$ be a Severi-Brauer variety of index $i(P)$. 
Then $\sym^{d}(P)$ is stably rational iff $i(P)\mid d$.
\end{cor}

Proof. If $d$ is not divisible by $i(P)$ then    $\sym^{d}(P)(k)=\emptyset$,  hence $\sym^{d}(P)$ is not stably rational.

To see the converse, assume that $i(P)\mid d$. By  (\ref{SB.basic.say}.5--6)
  $P$  is 
stably birational to a  Severi-Brauer variety
 $P'$   of dimension $d-1$.
Furthermore,  $\sym^{d}(P')$ is rational
by Theorem \ref{kra-sal.thm} and it is stably birational to
$\sym^{d}(P)$ by Corollary \ref{sym.of.prod.cor}. \qed

\begin{say}[Proof of Theorem \ref{sym.SB.thm}.1--6]
The easiest is (\ref{sym.SB.thm}.5). Given $d$ points in general position, they span a linear subspace of dimension $d-1$. This gives a natural map 
$\pi: \sym^d(P)\map \grass(\p^{d-1}, P)$. Let $K$ be the function field of
$\grass(\p^{d-1}, P)$ and $L_K\subset P_K$ the linear subspace corresponding to the generic point. Thus $L_K$ is a  
Severi-Brauer subvariety of dimension $d-1$.
The generic fiber of $\pi$ is $\sym^d(L_K)$ which is rational by
Theorem \ref{kra-sal.thm}.  Thus 
$\sym^d(P)\simb \grass(\p^{d-1}, P)\times \p^{d(d-1)}$.

Next we   show  (\ref{sym.SB.thm}.6) using the
 stable birational equivalences
$$
\begin{array}{rcl}
\sym^d(P)\times \sym^e(P)\times \p^{n(d,e)}&\sims&
\sym^d(P)\times \sym^e(P)\times  \sym^{(d,e)}(P)\\
&\sims&
\p^{nd}\times \p^{ne}\times  \sym^{(d,e)}(P).
\end{array}
$$
First let $K$ be the function field of $\sym^d(P)\times \sym^e(P)$.
Then $P_K$ has  0-cycles of degrees $d$ and $e$, thus 
is also has a 0-cycle of degree $(d, e)$. 
Thus  $\sym^{(d, e)}(P_K)$ is stably rational by Corollary \ref{kra-sal.cor},
proving the first part.

Similarly, let $L$ be the function field of $ \sym^{(d, e)}(P) $.
Then  $\sym^{d}(P_L)$  and $\sym^{e}(P_L)$ are stably
rational by Lemma \ref{kra-sal.cor}, proving the second part.

Using this for $e=i(P)$ gives that
$$
\sym^d(P)\times \sym^{i}(P)\sims \sym^{(d, i)}(P).
$$
Since $\sym^{i}(P)$ is stably rational by Corollary \ref{kra-sal.cor},
we (\ref{sym.SB.thm}.4).
Together with (\ref{sym.SB.thm}.5) this implies (\ref{sym.SB.thm}.2).

We have proved that every class in 
$\operatorname{MSym}(P)$ is stably birational to 
a symmetric power  $\sym^d(P)$ for some $ d\mid i(P)$.
Next we  show that these $\sym^d(P)$ are not stably birational
to each other.

Let $d<  e$  be different  divisors of $i(P)$. 
There is thus a prime $p$ such that $d=p^ad', e=p^ce'$ where
$a<c$ and $d', e'$ are not divisible by $p$. Let $p^b$ be the
largest $p$-power dividing $i(P)$.

By assumption $P$ has a $k'$ point for some field extension
$k'/k$ of degree $i(P)$.
Let $k''/k$ be the Galois closure of $k'/k$ and 
  $K$  the invariant subfield of a $p$-Sylow subgroup
 of $\gal(k''/k)$. 
Set $K'=k'K$.
Note that  $p$ does not divide $\deg (K/k)$ and $\deg (K'/K)=p^b$, hence
 $i(P_K)=p^b$.


Although $K'/K$ need not be Galois, 
the Galois group of its Galois closure is a $p$-group, hence nilpotent. Thus there is a
 subextension $K'\supset L\supset K$ of degree $p^{b-a}$. 
It is enough to show that   
 $\sym^d(P_L)$  and $\sym^e(P_L)$ are not stably birational
 over $L$.   
By (\ref{sym.SB.thm}.4),
$$
\sym^d(P_L)\sims \sym^{p^a}(P_L)\qtq{and} 
\sym^e(P_L)\sims \sym^{p^c}(P_L).
$$
Note that  $P_L$ has a point in $K'$ and $\deg (K'/L)=p^a$, hence  $i(P_L)=p^a$
and so $\sym^{p^a}(P_L) $ is stably rational by Corollary \ref{kra-sal.cor}.
By contrast
 $\sym^{p^c}(P_L) $ does not have  any $L$-points. Indeed, an $L$-point on
$\sym^{p^c}(P_L) $ would mean a 0-cycle of degree $p^c$ on $P_L$
hence a   0-cycle of degree $p^{b-a}p^c=p^{b-a+c}$ on $P_K$.
This is impossible since  $i(P_K)=p^b$ and $b-a+c<b$. 
Thus $\sym^{p^a}(P_L) $ and $\sym^{p^c}(P_L) $  are not stably birational.  
\qed
\end{say}

\begin{rem} It is possible that the stable birational equivalences in
Theorem \ref{sym.SB.thm} can be replaced by birational equivalences. 
For instance, it is possible that
$$
\sym^d(P)\simb \grass(\p^{(d, i)-1},P)\times \p^{m}\qtq{for suitable $m$.} 
$$
However, several steps in the proof naturally give only
stable birational equivalences and 
the difference between stable birational equivalence
and birational equivalence is not even understood for
Severi-Brauer varieties.
\end{rem}

We have used some general results on symmetric powers.

\begin{lem} \label{sym.of.prod}
Let $U$ be a positive dimensional, geometrically irreducible $k$-variety. Then
$$
\sym^m(U\times \p^r)\simb \sym^m(U)\times \p^{rm}.
$$
\end{lem}

Proof. There is a  natural projection map $ \sym^m(U\times \p^r)\to \sym^m(U)$. We claim that its generic fiber $F_{gen}$ is rational.
To construct it, set 
$$
L:=k\bigl(\sym^m(U)\bigr)\qtq{and}
K:=k\bigl(\sym^{m-1}(U)\times U\bigr).
$$
Here we think of $\sym^m(U)$ as $U^m/S_m$ and 
$\sym^{m-1}(U)\times U$ as $U^m/S_{m-1}$ where  $S_{m-1}\subset S_m$
are the permutations that fix the last factor. Thus $K/L$ is a
degree $m$ field extension and $F_{gen}\simb\Re_{K/L}(\p^r)$, the Weil restriction of $\p^r$ from $K$ to $L$.
Thus $F_{gen}$ is rational. \qed

\begin{cor} \label{sym.of.prod.cor} Let $U, V$ be 
positive dimensional, geometrically irreducible $k$-varieties. 
If $U\sims V$ then $\sym^m(U)\sims \sym^m (V)$ for every $m$.\qed
\end{cor}

As a consequence we see that $\sym^m(\p^r)$ is stably rational. In fact it is
rational; see \cite{MR0225774} for a very short proof.

\section{Moduli of  Severi-Brauer subvarieties}

We need some results on twisted line bundles and
maps between Severi-Brauer varieties.

\begin{defn}[Twisted line bundles]\label{twist.l.defn}
 Let $X$ be a geometrically normal, proper  $k$-variety.
A {\it twisted line bundle} of $X$ is  a line bundle $L$
on $X_{k^s}$ such that 
$L^{\sigma}\cong L$ for every $\sigma\in \gal(k^s/k)$. 
Equivalently, its class $[L]$ is a $k$-point of  $\pic(X)$.
For example, if $P$ is a Severi-Brauer variety then $\o_P(r)$ is a
twisted line bundle for every $r$.

Let $|L|$ denote the
irreducible component of the Hilbert scheme  (or Chow variety) of $X$
parametrizing subschemes $H\subset X$ such that $H_{k^s}$ is in the linear system $|L_{k^s}|$.
(See \cite[Chap.I.]{rc-book} for the Hilbert scheme  or the Chow variety.)
 This is clearly a Severi-Brauer variety.
There is a natural map $\iota_L:X\map |L|^{\vee}$ given by
$x\mapsto \{H: H\ni x\}$. 

Using this we define the {\it dual}  of  a Severi-Brauer variety $P$ as  $P^{\vee}:=|\o_P(1)|$.

 Let $\phi:X\map Y$ be a  map between geometrically normal, proper
 varieties and $L_Y$ a 
twisted line bundle on $Y$. Assume that either $\phi$ is a morphism or
$X$ is smooth. Then $\phi^*L_Y$ is a twisted line bundle on $X$ and
$|\phi^*L_Y|\sim |L_Y|$. 

Let $X,Y$ be  geometrically normal, proper
 varieties and $L_X, L_Y$ twisted line bundles on them. Let 
$\maps\bigl((X, L_X), (Y, L_Y)\bigr)$ denote the moduli space of
all maps $\phi:X\map Y$ such that $\phi^*L_Y\cong L_X$. 

If $P, Q$ are  Severi-Brauer varieties then we write
$$
\maps_d(Q, P):=\maps\bigl((Q, \o_Q(d)), (P, \o_P(1))\bigr).
$$
 Composing with $\iota_L$ gives an isomorphism 
$$
\maps\bigl((X, L), (P, \o_P(1))\bigr)\cong\maps_1\bigl(|L|^{\vee}, P\bigr).
$$
\end{defn}

\begin{say}[Severi-Brauer varieties II]\label{SB.B.say.2}
Let $P,Q$ be Severi-Brauer varieties. 
\begin{enumerate}
\item Their product     is defined as
$|\o_{P^{\vee}\times Q^{\vee}}(1,1)|\cong |\o_{P\times Q}(1,1)|^{\vee}$. 
I denote this by $P\otimes Q$. 
It is better to think of this as defined on Brauer-equivalence classes. This makes the set of
Brauer-equivalence classes into a group with identity $\p^0\sim \p^m$ and
inverse $P^{\vee}$.  The group is torsion, more precisely 
$P^{\otimes i(P)}\sim \p^0$.
(Frequently a smaller power of $P$ is trivial, the smallest such exponent is the {\it period.}) 

The group defined above  is isomorphic to  the {\it Brauer group} of $k$. 
(We will not use its cohomological description; see \cite{MR2266528}.)
\item $\maps_1(Q,P)\sim Q^{\vee}\otimes P$.  The natural map is
$\phi\mapsto\{(x,H): \phi(x)\in H\}\in |\o_{Q\times P^{\vee}}(1,1)|$.
\item   If $|L|$ is non-empty then $|L^m|\sim |L|^{\otimes m}$; this comes from
identifying the symmetric power of a vector space $V$ with the
subspace of symmetric tensors in $V^{\otimes m}$.  
\item $\maps_d(Q,P)\sim (Q^{\vee})^{\otimes d}\otimes P$; this follows from the previous two claims.
\item Again combining the previous two claims with (\ref{SB.basic.say}.1) we conclude that  there is a rational map  $P\map Q$ iff $Q$ is similar to
$P^{\otimes m}$ for some $m$.
(This is called Amitsur's theorem.)
\end{enumerate}
\end{say}

We next define the spaces of Severi-Brauer subvarieties of a 
Severi-Brauer variety. 
That is, given a Severi-Brauer variety $P$ we look at the subset of
the Chow variety $\chow(P)$ parametrizing subvarieties
$X\subset P$ whose normalization $\bar X$ is a Severi-Brauer variety.
For technical reasons it is better to work with   $\bar X\to P$. 

\begin{defn}
Fix integers $0\leq m\leq n$,   $1\leq d$ and a Severi-Brauer variety
$P$ of dimension $n$.  Let
$M^{\circ}_{\p^m}(P,d)$ denote the moduli space parametrizing morphisms 
$\phi:Q\to P$ satisfying the following assumptions.
\begin{enumerate}
\item $Q$ is a  Severi-Brauer variety of dimension $m$.
\item $\phi^*\o_P(1)\cong \o_Q(d)$.
\item Either $m<n$ and $\phi:Q\to \phi(Q)$ is birational or
$m=n$ and  every automorphism of the triple
$(\phi:Q\to P)$ that is the identity on $P$ is also the identity on $Q$.
\item Two such morphisms 
$\phi_i:Q_i\to P$ are identified if there is an isomorphism
$\tau:Q_1\cong Q_2$ such that $\phi_1=\phi_2\circ \tau$.
\end{enumerate}
The spaces $M^{\circ}_{\p^m}(P,d)$ are quasi-projective. They should be thought of as open subschemes of   the projective moduli spaces  of stable maps 
$\bar M_{\p^m}(P,d)$  \cite{ale-pairs}.
Since we are interested in their birational properties, these compactifications are not important to us.

(Comment on the notation. The moduli space of maps from $X$ to $Y$ is frequently denoted by $\maps(X, Y)$. However, for moduli of stable maps from a genus $g$ curve to $Y$ the usual notation is $M_g(Y,\beta)$ where $\beta$ is the homology class of the image.   If $Y=\p^n$ then $\beta$ is usually replaced by $\deg \beta$. Thus
$M^{\circ}_{\p^m}(P,d)$  follows mostly the
stable maps convention, except that the degree of
$\phi(Q)$ is $d^m$.)

Note that if $\phi:Q\to \phi(Q)$ is birational then every automorphism of
$\phi:Q\to P$ that is the identity on $P$ is also the identity on $Q$. This is why the most naive way of identifying two maps is adequate in (4). 
(As we discussed earlier, failure of this is one of the problems with $\bar M_0(P, 2)$ if $\dim P=1$.)

If $d=1$ then we get $M^{\circ}_{\p^m}(P,1)=\grass(\p^m, P)$ and if $m=1$
then the $M^{\circ}_{\p^1}(P,d)$ are open subschemes of the space of genus 0 stable maps  $\bar M_0(P,d)$.

These moduli spaces are closely related to the spaces of maps from
Definition \ref{twist.l.defn}:
$$
M^{\circ}_{\p^m}(P,d)\simb \maps_d(\p^m, P)/\aut(\p^m).
$$
The resulting map $\Pi:\maps_d(\p^m, P)\map M^{\circ}_{\p^m}(P,d)$ is not a product, not even birationally.   Indeed the fiber 
of $\Pi$ over a given $\phi:Q\to P$ is the space of isomorphisms
$\isom(\p^m, Q)$. This is a principal homogeneous space under
$\aut(\p^m)$ but it is not isomorphic to $\aut(\p^m)$ unless $Q$ is trivial.
\end{defn}

Our aim is to understand the spaces $M^{\circ}_{\p^m}(P,d)$
for arbitrary ground fields. This is  achieved only for $m=1$
but we have the following general periodicity property.

\begin{thm}\label{maps.perodic.thm}
 Let $P\sim P'$ be Brauer equivalent  Severi-Brauer varieties of dimensions
 $n, n'$.
Fix $ 0\leq m\leq \min\{n,n'\}$ and $1\leq d,d'$. Assume that
$d\equiv d'\mod (m+1)$. Then
$$
 M^{\circ}_{\p^m}(P,d)\sims M^{\circ}_{\p^m}(P', d').
$$
\end{thm}

Proof. The idea is  similar to the ``no-name method''  explained in 
\cite[Sec.4]{MR927970}, where it is attributed to Bogomolov and Lenstra. 

 Let  $ \isom_{\p^m}(d,P, d', P')$ denote the scheme parametrizing
triples
$$
\{(\phi:Q\to P); (\phi':Q'\to P'); \tau\}
$$
where 
 $(\phi:Q\to P)\in M^{\circ}_{\p^m}(P,d)$,
 $(\phi':Q'\to P')\in M^{\circ}_{\p^m}(P', d')$, and
 $\tau:Q\to  Q'$ is an isomorphism.
(No further assumptions on $\phi$ and $\phi'\circ \tau$.)
 We prove that
$$
M^{\circ}_{\p^m}(P,d)\sims \isom_{\p^m}(d,P,d', P')\sims M^{\circ}_{\p^m}(P', d'),
$$
using the  natural projections  
$$
\pi:\isom_{\p^m}(d,P,d', P')\to M^{\circ}_{\p^m}(P,d)\qtq{and}
\pi':\isom_{\p^m}(d,P,d', P')\to M^{\circ}_{\p^m}(P', d').
$$
 It is sufficient to show that their generic fibers are  rational. The roles of $d, d'$ are symmetrical,
thus it is enough to consider 
$\pi:\isom_{\p^m}(d,P,d',P')\to M^{\circ}_{\p^m}(P,d)$.

Note that the fiber of $\pi$ over  $(\phi:Q\to P)$
consists of pairs
$$
\{ (\phi':Q'\to P'); \tau\}
$$
where 
  $(\phi':Q'\to P')\in M^{\circ}_{\p^m}(P', d')$ and
 $\tau:Q\to  Q'$ is an isomorphism.
Specifying such a pair is the same as giving
$(\phi\circ \tau:Q\to P')\in \maps_d(Q,P')$. 
Thus the fiber of $\pi$ over  $(\phi:Q\to P)$ is isomorphic to
$\maps_{d'}(Q,P')$. 

Let $K$ be the function field of $M^{\circ}_{\p^m}(P,d)$. We thus have a
morphism  
$$
\phi_K:Q_K\to P_K \qtq{such that} \phi_K^*\o_{P_K}(1)\cong \o_{Q_K}(d).
$$
 By (\ref{SB.B.say.2}.4)
 $\phi_K$ corresponds to  a $K$-point of  
$ \maps_{d}(Q_K,P_K)\sim (Q_K^{\vee})^{\otimes d}\otimes P_K$.
Thus $(Q_K^{\vee})^{\otimes d}\otimes P_K$ is rational by (\ref{SB.basic.say}.1).
Furthermore, since $d\equiv d'\mod (m+1)$, we know that
$Q_K^{\otimes d'}\sim Q_K^{\otimes d}$ by (\ref{SB.B.say.2}.1), hence
$(Q_K^{\vee})^{\otimes d'}\otimes P'_K$ is stably rational by (\ref{SB.B.say.2}.1),
hence in fact rational by (\ref{SB.basic.say}.1). Therefore
$$
\maps_{d'}(Q_K,P'_K)\sim (Q_K^{\vee})^{\otimes d'}\otimes P'_K
$$
is also rational by (\ref{SB.basic.say}.1). 
 \qed

\begin{rem} There are a few more cases when one can get stable birational equivalences. For example, assume that  
 $d,d'$ and $\operatorname{per}(P)=\operatorname{per}(P')$, the period of $P$, are all relatively prime to $m+1$. Then
$$
 \phi_K^*\o_{P_K}\bigl( \operatorname{per}(P)\bigr)
\cong \o_{Q_K}\bigl( d'\cdot \operatorname{per}(P)\bigr)
$$ implies that $Q_K$ is trivial. 
Using this observation for $d'=1$ we obtain that
$$
 M^{\circ}_{\p^m}(P,d)\sims  M^{\circ}_{\p^m}( P',1)\cong \grass(\p^m, P')
\qtq{if}  \bigl(m+1, d\cdot   \operatorname{per}(P)\bigr)=1.
$$
\end{rem}

\medskip

As a consequence of Theorem \ref{maps.perodic.thm},
 in order to describe the 
stable birational types of   $M^{\circ}_{\p^m}(P,d)$,
it is sufficient to understand  $M^{\circ}_{\p^m}(P,d)$ for $d\leq m+1$.
There are two cases for which the answer is easy to derive.

\begin{lem} \label{0.1.lem}
 Let $P$ be a Severi-Brauer variety. Then
\begin{enumerate}
\item $M^{\circ}_{\p^m}(P,d)\sims\grass(\p^m, P) $ if $d\equiv 1\mod (m+1)$.
\item $ M^{\circ}_{\p^m}(P,d)\sims P$ if $d\equiv 0\mod (m+1)$ and
$(m+1)\mid 420$.
\end{enumerate}
\end{lem}

Proof. If  $d\equiv 1\mod (m+1)$ then
$M^{\circ}_{\p^m}(P,d)\sims M^{\circ}_{\p^m}( P,1)$ 
by Theorem \ref{maps.perodic.thm}
and,  essentially by definition, 
  $M^{\circ}_{\p^m}( P,1)=\grass(\p^m, P)$.

For the  second claim we check  the stable birational isomorphisms
$$
M^{\circ}_{\p^m}(P, m+1)\sims 
M^{\circ}_{\p^m}(P, m+1)\times P\sims P.
$$
First let $K$ be the function field of $M^{\circ}_{\p^m}(P, m+1) $.
We need to show that $P_K$ is trivial.
By assumption, there is a $K$-map in
$\maps_{m+1}(Q_K, P_K)$ where $\dim Q_K=m$. By (\ref{SB.B.say.2}.4)
this corresponds to a $K$-point of   $(Q_K^{\vee})^{\otimes m+1}\otimes P_K$. 
By (\ref{SB.B.say.2}.1)  $Q_K^{\otimes m+1}$ is trivial, so
$\maps_{m+1}(Q_K, P_K)\sim P_K$ and so $P_K$ is trivial.

For the second part, let $L$ be the function field of $P$.
Then $P_L$ is trivial, hence
$$
M^{\circ}_{\p^m}(P_L, m+1)\cong M^{\circ}_{\p^m}(\p^n_L, m+1)\simb
 \p\bigl(H^0(\p^m_L, \o_{\p^m}(m+1))^{n+1}\bigr)/\PGL_{m+1}.
$$
It is conjectured that this quotient if always stably 
rational, but this seems to be known only when $(m+1)\mid 420$; see 
\cite[p.316]{MR1958908} and the references there. \qed


\medskip

We give a  geometric proof that 
the space of conics $\bar M_0(P,2)$ is stably
birational to $P$ for $\dim P\geq 2$.

\begin{say}[Conics in Severi-Brauer varieties]\label{lin.con.say}
 We compute, in 2 different ways,  the space $T$ 
parametrizing  triples
$(C,\ell_1, \ell_2)$ where the $\ell_i$ are secant lines of $C$. 

Forgetting the lines gives a map to $\bar M_0(P,2)$. Let $C_K$ be the conic corresponding to the generic point of $\bar M_0(P,2)$. A secant line of $C_K$ is
determined by $\sym^2C_K\cong \p^2_K$. Thus
$T\simb \bar M_0(P,2)\times \p^4$.

The secants lines $\ell_i$ meet at a unique point; this gives a map
$T\to P$.  Given $p\in P$, the fiber is obtained by first picking
2 points in $\p(T_pP)\cong \p^{n-1}_{k(p)}$. Once we have 2 lines, they determine a plane $\langle\ell_1, \ell_2\rangle$ and the  5-dimensional linear system
$|\ell_1+ \ell_2|$ on  the plane $\langle\ell_1, \ell_2\rangle$ gives 
the conics. Thus 
$T\simb P\times \p^{2n+3}$ and hence  $\bar M_0(P,2)\sims P$.
\end{say}

\begin{say}[Proof of Theorem \ref{sym.SB.thm}.7--8]
If $d=2e$ is even then
$$
\bar M_0(P,2e)\simb M_{\p^1}(P,2e)\sims M_{\p^1}(P,2)\simb \bar M_0(P,2),
$$
where the birational equivalences are by definition and
the stable birational equivalence holds by Theorem \ref{maps.perodic.thm}.
Next $\bar M_0(P,2)\sims P$ follows either from Lemma
\ref{0.1.lem}.2 or from Paragraph \ref{lin.con.say}.
This gives (\ref{sym.SB.thm}.7). 

Similarly, if $d=2e+1$ is odd then (\ref{sym.SB.thm}.8) follows from 
$$
\bar M_0(P, 2e+1)\simb M_{\p^1}(P,2e+1)\sims M_{\p^1}( P,1)=\grass(\p^1, P). \qed
$$
\end{say}

\begin{rem}
So far  we have  worked with a fixed Severi-Brauer variety $P$, but it would be interesting to understand how the $\operatorname{MSym}(P)$ for different  Severi-Brauer varieties interact with each other.  

For example, assume that $P, Q$ are  Severi-Brauer varieties
such that $\ind(P)$ and $\ind(Q)$ are relatively prime.
We claim that $\sym^d(P)$ and $\sym^e(Q)$ are stably birational
to each other iff they are both stably rational. 

To see this, assume that $\sym^e(Q)$ is not stably rational.
By 
Corollary \ref{kra-sal.cor} this holds iff $\ind(Q)\nmid e$.
Set $K:=k(P)$. Since $\ind(P)$ and $\ind(Q)$ are relatively prime,
$\ind(Q_K)=\ind(Q)$, thus $\sym^e(Q_K)$ is not stably rational by
Corollary \ref{kra-sal.cor}. By contrast, $P_K$ is trivial hence
$\sym^d(P_K)$ is stably rational; even rational by \cite{MR0225774}. Thus
$\sym^d(P)$ and $\sym^e(Q)$ are not stably birational
to each other.

Further steps in this direction are in
Theorem \ref{maps.perodic.thm}; see also \cite{k-conic, MR2439423} for
related questions.
\end{rem}

\begin{ack} I thank I.~Coskun, P.~Gille, D.~Krashen, M.~Lieblich, 
T.~Szamuely  for helpful
comments and H.-C.~von~Bothmer for calling my attention to \cite{MR1958908}.
Partial financial support    was provided  by  the NSF under grant number
 DMS-1362960.
\end{ack}


\begin{thebibliography}{Mat68}

\bibitem[Ale96]{ale-pairs}
Valery Alexeev, \emph{Moduli spaces {$M\sb {g,n}(W)$} for surfaces},
  Higher-dimensional complex varieties ({T}rento, 1994), de Gruyter, Berlin,
  1996, pp.~1--22.

\bibitem[Dol87]{MR927970}
Igor~V. Dolgachev, \emph{Rationality of fields of invariants}, Algebraic
  geometry, {B}owdoin, 1985 ({B}runswick, {M}aine, 1985), Proc. Sympos. Pure
  Math., vol.~46, Amer. Math. Soc., Providence, RI, 1987, pp.~3--16.
  \MR{927970}

\bibitem[For02]{MR1958908}
Edward Formanek, \emph{The ring of generic matrices}, J. Algebra \textbf{258}
  (2002), no.~1, 310--320, Special issue in celebration of Claudio Procesi's
  60th birthday. \MR{1958908}

\bibitem[GS06]{MR2266528}
Philippe Gille and Tam{\'a}s Szamuely, \emph{Central simple algebras and
  {G}alois cohomology}, Cambridge Studies in Advanced Mathematics, vol. 101,
  Cambridge University Press, Cambridge, 2006. \MR{2266528}

\bibitem[Hog09]{MR2439423}
Amit Hogadi, \emph{Products of {B}rauer-{S}everi surfaces}, Proc. Amer. Math.
  Soc. \textbf{137} (2009), no.~1, 45--50. \MR{2439423}

\bibitem[Kol96]{rc-book}
J{\'a}nos Koll{\'a}r, \emph{Rational curves on algebraic varieties}, Ergebnisse
  der Mathematik und ihrer Grenzgebiete. 3. Folge., vol.~32, Springer-Verlag,
  Berlin, 1996.

\bibitem[Kol05]{k-conic}
\bysame, \emph{Conics in the {G}rothendieck ring}, Adv. Math. \textbf{198}
  (2005), no.~1, 27--35. \MR{2183248 (2006k:14064)}

\bibitem[KS04]{MR2090670}
Daniel Krashen and David~J. Saltman, \emph{Severi-{B}rauer varieties and
  symmetric powers}, Algebraic transformation groups and algebraic varieties,
  Encyclopaedia Math. Sci., vol. 132, Springer, Berlin, 2004, pp.~59--70.
  \MR{2090670 (2005k:14024)}

\bibitem[Mat68]{MR0225774}
Arthur Mattuck, \emph{The field of multisymmetric functions}, Proc. Amer. Math.
  Soc. \textbf{19} (1968), 764--765. \MR{0225774}

\end{thebibliography}
\def\cprime{$'$} \def\cprime{$'$} \def\cprime{$'$} \def\cprime{$'$}
  \def\cprime{$'$} \def\cprime{$'$} \def\dbar{\leavevmode\hbox to
  0pt{\hskip.2ex \accent"16\hss}d} \def\cprime{$'$} \def\cprime{$'$}
  \def\polhk#1{\setbox0=\hbox{#1}{\ooalign{\hidewidth
  \lower1.5ex\hbox{`}\hidewidth\crcr\unhbox0}}} \def\cprime{$'$}
  \def\cprime{$'$} \def\cprime{$'$} \def\cprime{$'$}
  \def\polhk#1{\setbox0=\hbox{#1}{\ooalign{\hidewidth
  \lower1.5ex\hbox{`}\hidewidth\crcr\unhbox0}}} \def\cdprime{$''$}
  \def\cprime{$'$} \def\cprime{$'$} \def\cprime{$'$} \def\cprime{$'$}
\providecommand{\bysame}{\leavevmode\hbox to3em{\hrulefill}\thinspace}
\providecommand{\MR}{\relax\ifhmode\unskip\space\fi MR }
\providecommand{\MRhref}[2]{%
  \href{http://www.ams.org/mathscinet-getitem?mr=#1}{#2}
}
\providecommand{\href}[2]{#2}


\noindent  Princeton University, Princeton NJ 08544-1000

{\begin{verbatim} kollar@math.princeton.edu\end{verbatim}}

\end{document}